\documentclass[11pt]{article}
\usepackage[a4paper,margin=2.5cm]{geometry}
\usepackage{indentfirst}
\usepackage{amsmath,amssymb,amsthm,mathtools}
\usepackage{enumitem,microtype}
\usepackage{xcolor,needspace}
\usepackage[
	colorlinks=true,
	allcolors=black,
	bookmarksnumbered=true,
	pagebackref=true
]{hyperref}
\usepackage[normalem]{ulem}

\renewcommand*{\backref}[1]{}
\renewcommand*{\backrefalt}[4]{%
	\ifcase#1\relax
	\or
		(Cited on page~#2.)%
	\else
		(Cited on pages~#2.)%
	\fi
}

\DeclareMathOperator{\soc}{soc}
\DeclareMathOperator{\rad}{rad}
\DeclareMathOperator{\Hom}{Hom}

\DeclareMathOperator{\Aut}{Aut}
\newcommand{\defn}[1]{\textbf{#1}}

\newtheorem{theorem}{Theorem}[section]
\newtheorem{proposition}[theorem]{Proposition}
\newtheorem{lemma}[theorem]{Lemma}
\newtheorem{corollary}[theorem]{Corollary}
\theoremstyle{definition}
\newtheorem{definition}[theorem]{Definition}
\newtheorem{remark}[theorem]{Remark}

\setlist[enumerate]{label=\textup{(\arabic*)},itemsep=3pt,topsep=4pt}
\title{A characterisation of almost gentle algebras}
\author{Pengyun Chen$^{a}$ and Bohan Xing$^{a,*}$}
\date{}
\hypersetup{
	pdftitle={A characterisation of almost gentle algebras},
	pdfauthor={Pengyun Chen and Bohan Xing}
}

\begin{document}
	\maketitle

	\begingroup
		\renewcommand{\thefootnote}{}
		\makeatletter
		\renewcommand{\@makefntext}[1]{%
			\setlength{\parindent}{2em}%
			\noindent\hspace*{\parindent}#1}
		\makeatother
		\begin{NoHyper}
			\footnotetext{%
				The version of \today.
				\par
				$^{a}$Pengyun Chen and Bohan Xing,
				School of Mathematical Sciences,
				Laboratory of Mathematics and Complex Systems,
				Beijing Normal University,
				Beijing 100875, P. R. China.
				\par
				E-mail addresses:
				\texttt{pychen@mail.bnu.edu.cn} (P. Chen);
				\texttt{bhxing@mail.bnu.edu.cn} (B. Xing).
				\par
				$^{*}$Corresponding author.
			}
		\end{NoHyper}
	\endgroup

	\begin{abstract}
		A finite-dimensional algebra over an algebraically closed field is gentle if and only if its trivial extension is special biserial.
		We prove the corresponding characterisation for almost gentle algebras, that is, an algebra over such a field is almost gentle if and only if its trivial extension is special multiserial, which answers a question of Green and Schroll.
		As a consequence, the property of having a special multiserial trivial extension is invariant under stable equivalence.
	\end{abstract}

\section{Introduction}
\label{sec:introduction}

	It is well known that gentle algebras can be characterised by their trivial extensions.
	To be more explicit, over an algebraically closed field, a finite-dimensional algebra $A$ is gentle if and only if $T(A)$ is special biserial~\cite{PogorzalySkowronskiBiserial}.
	Here $T(A)=A\ltimes D(A)$ is the trivial extension of $A$.
	The relation between these classes can also be described by admissible cuts: every Brauer graph algebra with multiplicities one is the trivial extension of a gentle algebra obtained by an admissible cut~\cite{SchrollTrivial}.

	Green and Schroll generalised gentle algebras to almost gentle algebras~\cite{GreenSchrollAlmost}, which allow arbitrarily many arrows to enter or leave a vertex and impose no restriction on the number of zero continuations of an arrow on either side.
	Note that the relations remain quadratic monomial, and each arrow still has at most one nonzero continuation on either side.
	Thus the special biserial condition is relaxed to the special multiserial condition~\cite{vonHohneWaschbusch,GreenSchrollMultiserial}.
	These classes of algebras also provide positive cases of classical homological conjectures; see~\cite{Sch19,WHL25,CX26}.
	Green and Schroll proved that the trivial extension of an almost gentle algebra is special multiserial, and asked whether the converse holds~\cite[Corollary~4.4 and Question~4.5]{GreenSchrollAlmost}.
	We answer their question as follows.

	\begin{theorem}
		\label{main:trivial-extension}
		\label{thm:almost-gentle-trivial-special-multiserial}
		Let $A$ be a finite-dimensional algebra over an algebraically closed field.
		Then $A$ is almost gentle if and only if $T(A)$ is special multiserial.
	\end{theorem}

	Combining Theorem~\ref{main:trivial-extension} with our results~\cite{CX26}, we obtain consequences for stable equivalences.
	The property that $T(A)$ is special multiserial is preserved under stable equivalence of $A$ (Corollary~\ref{cor:stable-special-multiserial}).
	Moreover, if basic algebras $A$ and $B$ without semisimple direct summands are stably equivalent of Morita type and $T(A)$ is special multiserial, then $A\cong B$ and $T(A)\cong T(B)$ (Corollary~\ref{cor:morita-type-trivial-extensions}).
	This gives a positive answer in this case to the question of whether taking trivial extensions preserves stable equivalences of Morita type.
	Liu, Zhou and Zimmermann~\cite{LZZ17} showed that the answer is negative in general.

	The paper is organised as follows.
	In Section~\ref{sec:preliminaries}, we recall the definitions and the results on trivial extensions and admissible cuts used in the proof.
	In Section~\ref{sec:trivial-extensions}, we prove Theorem~\ref{main:trivial-extension}.
	Section~\ref{sec:stable-equivalences} gives the consequences for stable equivalences.

\section{Preliminaries}
\label{sec:preliminaries}

	Throughout, $\Bbbk$ is an algebraically closed field and $\Bbbk^\times$ denotes its multiplicative group.
	All algebras are finite-dimensional associative $\Bbbk$-algebras with identity, and modules are left modules unless stated otherwise.
	We write $D=\Hom_\Bbbk(-,\Bbbk)$ for the usual duality and $\rad A$ and $\soc A$ for the radical and socle of an algebra $A$.
	For a finite quiver $Q$, let $Q_0$ and $Q_1$ denote its vertex and arrow sets, and let $s,t:Q_1\to Q_0$ be the source and target maps.
	We compose paths from right to left.
	An ideal $I$ of $\Bbbk Q$ is admissible if $R^n_Q\subseteq I\subseteq R_Q^2$ for some $n\geq2$, where $R_Q$ is the arrow ideal. 
	We use the same notation for a path and its class in $\Bbbk Q/I$, and write $e_i$ for the trivial path at $i\in Q_0$.

	\subsection{Special multiserial and almost gentle algebras}
	
		\begin{definition}[{\cite{vonHohneWaschbusch,GreenSchrollMultiserial,GreenSchrollAlmost}}]
			\label{def:special-multiserial}
			\label{def:almost-gentle}
			An algebra is \defn{special multiserial} if it is Morita equivalent to a bound quiver algebra $\Bbbk Q/I$ satisfying
			\begin{enumerate}
				\item[$(M)$] For every arrow $\alpha\in Q_1$, there is at most one arrow $\beta$ with $s(\beta)=t(\alpha)$ and $\beta\alpha\notin I$, and at most one arrow $\gamma$ with $t(\gamma)=s(\alpha)$ and $\alpha\gamma\notin I$.
			\end{enumerate}
			Such a presentation is called a \defn{special multiserial presentation}.
			An algebra is \defn{almost gentle} if it admits, up to Morita equivalence, a special multiserial presentation in which $I$ is generated by paths of length two.
		\end{definition}

		An algebra is \defn{monomial} if it has a presentation $\Bbbk Q/I$ with $I$ generated by paths.
		For such a presentation, the paths not in $I$ form a basis of $\Bbbk Q/I$.
		A nonzero path is \defn{maximal} if it cannot be extended on either side by an arrow to a nonzero path.

	\subsection{Trivial extensions and admissible cuts}
	\label{subsec:cycles}

		For an algebra $\Gamma$ and a $\Gamma$-$\Gamma$-bimodule $M$, the \defn{split square-zero extension} $\Gamma\ltimes M$ has underlying vector space $\Gamma\oplus M$ and multiplication
		\[
			(s,m)(s',m')=(ss',sm'+ms').
		\]
		In particular, the \defn{trivial extension} of $A$ is $T(A)=A\ltimes D(A)$, where $(afb)(x)=f(bxa)$ for $a,b,x\in A$ and $f\in D(A)$; see for example in~\cite[Section~2]{FernandezEtAlTrivialExtensions}.
		We view $A$ and $D(A)$ as subspaces of $T(A)$ via the embeddings $a\mapsto(a,0)$ and $f\mapsto(0,f)$, respectively.
		Then $D(A)$ is a square-zero ideal and $T(A)/D(A)\cong A$.

		Note that every trivial extension is symmetric: the linear map
		\[
			\tau:T(A)\longrightarrow\Bbbk,\qquad (a,f)\longmapsto f(1)
		\]
		satisfies $\tau(xy)=\tau(yx)$, and the form $(x,y)\mapsto\tau(xy)$ is nondegenerate.
		Trivial extensions respect finite products and Morita equivalence.
		Indeed, if $e\in A$ is a full idempotent, that is, $AeA=A$, then $e$ is also full in $T(A)$, and restriction of functionals gives
		\[
			eT(A)e\cong eAe\ltimes eD(A)e\cong T(eAe).
		\]

		Let $B$ be a basic connected symmetric special multiserial algebra with $\rad^2B\neq0$.
		We use the following \defn{cycle presentation} of $B$~\cite[Section 4]{GreenSchrollMultiserial}.
		In the Gabriel quiver $Q$ of $B$, choose oriented cycles $C$, without repeated arrows and up to cyclic permutation, whose arrow sets partition $Q_1$.
		Their cyclic permutations are the \defn{special cycles}.
		Each cycle $C$ has an associated positive integer $m_C$, called its \defn{multiplicity}, such that $C^{m_C}$ is a nonzero socle element.
		We write $|C|$ for the number of arrows in $C$, so that $C^{m_C}$ has length $m_C|C|$.
		We require $m_C|C|\geq2$.
		The nonzero paths of positive length are precisely the subpaths of the cyclic permutations of $C^{m_C}$.
		At each vertex $i$, the full powers based at $i$ are identified with an element $s_i$ spanning $e_i\soc B e_i$.
		For every vertex $i$ and every arrow $\alpha$, $s_i\alpha=\alpha s_i=0$.
		The vertex idempotents, the distinct proper subpaths of positive length, and the elements $s_i$ form a basis of $B$.

		An \defn{admissible cut} $\Delta$ contains exactly one arrow from each cycle up to cyclic permutation.
		We need the following result only when all multiplicities are one.

		\begin{theorem}[{\cite[Theorems~5.2 and~5.4]{GreenSchrollAlmost}}]
			\label{thm:cycle-admissible-cut}
			Let $B$ have a cycle presentation with all multiplicities one, and let $\Delta$ be an admissible cut.
			Then $\Lambda=B/\langle\Delta\rangle$ is almost gentle and $T(\Lambda)\cong B$.
		\end{theorem}

		Conversely, for a connected non-semisimple almost gentle algebra $\Lambda=\Bbbk Q/I$, a cycle presentation of $T(\Lambda)$ is obtained by adding an arrow $\beta_p:t(p)\to s(p)$ for each maximal nonzero path $p$ of $\Lambda$.
		The arrow $\beta_p$ is represented by the dual path functional $p^*\in D(\Lambda)$, and the cycles are the cyclic permutations of $\beta_p p$, all with multiplicity one~\cite[Theorem~4.3]{GreenSchrollAlmost}.

\section{Proof of the main theorem}
\label{sec:trivial-extensions}

	\subsection{Isomorphic trivial extensions}

		A split square-zero extension $\Lambda=\Gamma\ltimes M$ has the grading $\Lambda=\Lambda_0\oplus \Lambda_1$ with $\Lambda_0=\Gamma$ and $\Lambda_1=M$.
		The next lemma allows us to choose homogeneous arrow generators without changing a monomial relation ideal.

		\begin{lemma}
			\label{lem:monomial-grading-presentation}
			Let $\Lambda=\Bbbk Q/I$ be a monomial algebra with a grading $\Lambda=\Lambda_0\oplus \Lambda_1$ supported in degrees zero and one.
			Then there is an automorphism of $\Lambda$ sending each vertex into $\Lambda_0$ and each arrow to a homogeneous element.
		\end{lemma}

		\begin{proof}
			If $\Lambda_1=0$, there is nothing to prove. 
			Since $\Lambda_1$ is nilpotent and $\Lambda/\Lambda_1\cong \Lambda_0$, a complete set of primitive orthogonal idempotents of $\Lambda_0$ is also such a set in $\Lambda$. 
			Thus, after composing the presentation with an inner automorphism, we may assume that $e_i\in \Lambda_0$ for every $i\in Q_0$.
			Write $\Aut(\Lambda,e_i)$ for the group of $\Bbbk$-algebra automorphisms of $\Lambda$ fixing every $e_i$.
			This is a linear algebraic group.
			An algebraic subgroup isomorphic to $(\Bbbk^\times)^r$ is called a \defn{torus}; it is maximal if it is not properly contained in another torus.

			Since $I$ is monomial, each choice of nonzero scalars
			$\lambda=(\lambda_\alpha)_{\alpha\in Q_1}$ defines an automorphism
			\[
				d_\lambda:\Lambda\longrightarrow \Lambda,\qquad
				e_i\longmapsto e_i,\quad
				\alpha\longmapsto\lambda_\alpha\alpha.
			\]
			These automorphisms form a torus
			\[
				T=\{d_\lambda\mid
				\lambda\in(\Bbbk^\times)^{Q_1}\}
				\subseteq\Aut(\Lambda,e_i).
			\]
			We show that $T$ is maximal.
			Let $u$ belong to a torus containing $T$.
			Since a torus is isomorphic to $(\Bbbk^\times)^r$ for some integer $r\geq 0$, it is commutative, and hence for every arrow $\alpha$ we have
			\[
				d_\lambda(u(\alpha))
				=u(d_\lambda(\alpha))
				=\lambda_\alpha u(\alpha)
			\]
			for every choice of $\lambda$.

			Since $\Lambda$ is a monomial algebra, the paths outside $I$ form a basis of $\Lambda$.
			Note that on a path $p=\alpha_\ell\cdots\alpha_1$, the automorphism $d_\lambda$ acts by multiplication by $\lambda_{\alpha_\ell}\cdots\lambda_{\alpha_1}$.
			Since the scalars $\lambda_\beta$ vary independently, the only basis path multiplied by $\lambda_\alpha$ for every $\lambda$ is the arrow $\alpha$ itself.
			Comparing coefficients in the preceding equality therefore gives $u(\alpha)\in\Bbbk\alpha$.
			Thus $u$ fixes the vertices and rescales each arrow, so $u\in T$.
			Hence $T$ is maximal.

			For each $t\in\Bbbk^\times$, there is an automorphism
			\[
				h_t:\Lambda\longrightarrow\Lambda,\qquad
				x_0+x_1\longmapsto x_0+tx_1
				\quad(x_j\in\Lambda_j).
			\]
			They fix every $e_i$ and form a torus in $\Aut(\Lambda,e_i)$.
			By \cite[Theorem~17.10]{MilneAlgebraicGroups}, every torus is contained in a maximal torus, and maximal tori are conjugate over an algebraically closed field.
			Consequently, there is
			$\varphi\in\Aut(\Lambda,e_i)$ such that
			\[
				\varphi h_t\varphi^{-1}\in T
				\qquad\text{for every }t\in\Bbbk^\times.
			\]
			Since every element of $T$ sends each arrow to a scalar multiple of itself, this implies
			\[
				h_t\bigl(\varphi^{-1}(\alpha)\bigr)
				\in\Bbbk\varphi^{-1}(\alpha)
				\qquad\text{for every }t\in\Bbbk^\times.
			\]
			Write $\varphi^{-1}(\alpha)=x_0+x_1$ with $x_j\in \Lambda_j$.
			Then $x_0+tx_1\in\Bbbk(x_0+x_1)$ for every $t$.
			Taking $t\neq1$ and comparing the two homogeneous components shows that $x_0=0$ or $x_1=0$.
			Thus every $\varphi^{-1}(\alpha)$ is homogeneous.

			Finally, compose the presentation map $\Bbbk Q\to \Lambda$ with $\varphi^{-1}$.
			The vertex images remain in $\Lambda_0$, the arrow images are homogeneous, and the kernel remains $I$.
		\end{proof}

		\begin{theorem}[{\cite{WakamatsuTrivialExtensions}, see also \cite[Theorem~4.1]{FernandezEtAlTrivialExtensions}}]
			\label{thm:wakamatsu-trivial-extensions}
			Let $\Lambda$ and $\Lambda'$ be basic algebras.
			Then $T(\Lambda)\cong T(\Lambda')$ if and only if there are an algebra $\Gamma$ and a $\Gamma$-$\Gamma$-bimodule $M$ such that
			\[
				\Lambda\cong \Gamma\ltimes M,
				\qquad
				\Lambda'\cong \Gamma\ltimes D(M).
			\]
		\end{theorem}

		\begin{proposition}
			\label{prop:almost-gentle-trivial-invariance}
			Let $\Lambda$ and $\Lambda'$ be basic algebras with $T(\Lambda)\cong T(\Lambda')$.
			If $\Lambda$ is almost gentle, then $\Lambda'$ is almost gentle.
		\end{proposition}

		\begin{proof}
			We may assume that $\Lambda$ is connected.
			If $\Lambda=\Bbbk$, then $\dim_\Bbbk \Lambda'=1$ and there is nothing to prove.
			Otherwise, Theorem~\ref{thm:wakamatsu-trivial-extensions} gives
			\[
				\Lambda=\Gamma\ltimes M,\qquad \Lambda'\cong \Gamma\ltimes D(M).
			\]
			Give $\Lambda$ the grading $\Lambda_0=\Gamma$ and $\Lambda_1=M$, and fix an almost gentle presentation $\Lambda\cong\Bbbk Q/I$.
			By Lemma~\ref{lem:monomial-grading-presentation}, we may choose the presentation map so that the vertex idempotents have degree zero and the arrow generators are homogeneous, while keeping $Q$ and $I$ unchanged.
			Thus the resulting presentation is still almost gentle.

			Extend this grading to $T(\Lambda)$ by setting
			\[
				T(\Lambda)_0=\Gamma\oplus D(M),\qquad T(\Lambda)_1=M\oplus D(\Gamma),
			\]
			where $D(M)$ and $D(\Gamma)$ are embedded in $D(\Lambda)$ by extending functionals by zero.
			For a nonzero path $p$ in $\Lambda$, its dual basis element $p^*$ has degree $1-\deg p$.
			Thus, in the cycle presentation of $T(\Lambda)$ recalled in Subsection~\ref{subsec:cycles}, each returning arrow $\beta_p$ is homogeneous and every special cycle $\beta_p p$ has degree one.

			Let $\Delta$ be the set of arrows of degree one.
			Each special cycle contains exactly one such arrow, so $\Delta$ is an admissible cut.
			Moreover, $T(\Lambda)_1$ is spanned by paths of degree one, each of which contains an arrow of $\Delta$.
			Hence $\langle\Delta\rangle=T(\Lambda)_1$, and
			\[
				T(\Lambda)/\langle\Delta\rangle\cong T(\Lambda)_0
				\cong \Gamma\ltimes D(M)\cong \Lambda'.
			\]
			It follows from Theorem~\ref{thm:cycle-admissible-cut} that $\Lambda'$ is almost gentle.
		\end{proof}

	\subsection{Multiplicities in a cycle presentation}

		We now show that a special multiserial trivial extension admits a cycle presentation with all multiplicities one.

		\begin{proposition}
			\label{prop:multiplicity-one}
			Let $A$ be basic and connected, and suppose that $B=T(A)$ is special multiserial with $\rad^2B\neq0$.
			Then $B$ admits a cycle presentation with all multiplicities one.
		\end{proposition}

		\begin{proof}
			Choose a cycle presentation of $B$ as in Subsection~\ref{subsec:cycles}.
			After applying an inner automorphism of $B$ to the chosen presentation, we may assume that its vertex idempotents lie in $A$.
			Set
			\[
				K=0\oplus D(A),\qquad
				\pi:B\longrightarrow B,\qquad(a,f)\longmapsto(a,0).
			\]
			Then $\pi$ fixes the vertex idempotents, $\ker\pi=K$, $K^2=0$, and
			\[
				\soc B=0\oplus D(A/\rad A)\subseteq K\subseteq\rad B.
			\]

			We first show that every cycle of multiplicity greater than one is a loop of multiplicity two.
			Suppose that $C$ has length $r$ and multiplicity $m>1$, with $rm\geq3$.
			Fix an arrow $\alpha$ of $C$, and let $c_\alpha(d)$ denote the coefficient of $\alpha$ in $d\in K$ with respect to the basis associated with the chosen cycle presentation.

			If $r\geq2$, write a cyclic permutation of $C$ as $\alpha b$. The path $\alpha b\alpha$ has length $r+1<rm$, so it is a basis element that does not lie in $\soc B$. The coefficient of $\alpha b\alpha$ in $dbd$ is $c_\alpha(d)^2$, since terms involving longer paths have length greater than $r+1$ and condition $(M)$ excludes all other arrows.
			Since $dbd=0$, we obtain $c_\alpha(d)^2=0$.
			If $r=1$, then $m\geq3$, and the coefficient of $\alpha^2$ in $d^2=0$ gives the same conclusion.
			Hence $c_\alpha(d)=0$ for every $d\in K$ and every arrow $\alpha$ of $C$.

			Write $\pi(\alpha)=\alpha-d_\alpha$ with $d_\alpha\in K$ and expand $\pi(C^m)$.
			Terms containing at least two factors from $K$ vanish because $KBK=0$. Consider a term containing exactly one factor $d_\alpha$. By condition $(M)$, a path occurring in the expansion of $d_\alpha$ can give a nonzero product only if it is $\alpha$ or a longer path along $C$.
			The coefficient of $\alpha$ in $d_\alpha$ is zero, while a longer path makes the total length exceed $rm$ and hence gives zero.
			Socle terms are annihilated by the neighbouring arrows.
			Therefore
			\[
				\pi(C^m)=C^m\neq0,
			\]
			contrary to $C^m\in\soc B\subseteq\ker\pi$.

			It remains to consider loops $\alpha$ of multiplicity two at a vertex $i$, so that $\alpha^2=s_i$.
			If $\operatorname{char}\Bbbk=2$, the canonical symmetrizing trace satisfies
			\[
				\tau((a,f)^2)=2f(a)=0.
			\]
			But $\tau(s_i)\neq0$ by nondegeneracy, since $s_iB=\Bbbk s_i$.
			Thus no such loop occurs in characteristic two.

			Assume that $\operatorname{char}\Bbbk\neq2$.
			Fix a vertex $i$, and let $q_i$ be the number of multiplicity-two loops at $i$. For each cycle $C$ of multiplicity one, let $t_i(C)$ be the number of arrows of $C$ starting at $i$. Each ordered pair of distinct such arrows determines a proper nonzero path of positive length along $C$ from $i$ to $i$, so there are $t_i(C)(t_i(C)-1)$ such paths. Counting the basis elements in $e_iBe_i$, we obtain
			\[
				2\dim_\Bbbk e_iAe_i
				=\dim_\Bbbk e_iBe_i
				=2+q_i+\sum_C t_i(C)\bigl(t_i(C)-1\bigr),
			\]
			where the sum runs over the multiplicity-one cycles.
			Thus $q_i$ is even.

			Pair the multiplicity-two loops at each vertex.
			For a pair $\alpha,\beta$ at $i$, we have $\alpha^2=\beta^2=s_i$ and $\alpha\beta=\beta\alpha=0$.
			Choose $\iota\in\Bbbk$ with $\iota^2=-1$ and put
			\[
				u=\alpha+\iota\beta,\qquad
				v=\frac{\alpha-\iota\beta}{2}.
			\]
			The linear change from $(\alpha,\beta)$ to $(u,v)$ is invertible, and
			\[
				u^2=v^2=0,\qquad uv=vu=s_i.
			\]
			Both $u$ and $v$ annihilate every unchanged arrow on either side. Thus the pair $\alpha,\beta$ is replaced by the length-two cycle $uv$ of multiplicity one.
			Repeating this for every pair gives the required presentation.
		\end{proof}

		\smallskip
		\smallskip
		\begin{proof}[Proof of Theorem~\ref{main:trivial-extension}]
			Trivial extensions respect Morita equivalence and finite products, so we may assume that $A$ is basic and connected.

			\smallskip
			\noindent\underline{$(\Rightarrow)$}\quad
			If $A$ is almost gentle, then $T(A)$ is special multiserial by \cite[Corollary~4.4]{GreenSchrollAlmost}.

			\smallskip
			\noindent\underline{$(\Leftarrow)$}\quad
			Suppose that $B=T(A)$ is special multiserial. If $\rad^2B=0$, then $\rad A=0$. Hence $A\cong\Bbbk$ in this case and is almost gentle.
			If $\rad^2B\neq0$, Proposition~\ref{prop:multiplicity-one} gives a cycle presentation of $B$ with all multiplicities one.
			Choose an admissible cut $\Delta$ and put $\Lambda=B/\langle\Delta\rangle$.
			By Theorem~\ref{thm:cycle-admissible-cut}, the algebra $\Lambda$ is almost gentle and $T(\Lambda)\cong B=T(A)$.
			Proposition~\ref{prop:almost-gentle-trivial-invariance} now shows that $A$ is almost gentle.
		\end{proof}

\section{Consequences for stable equivalences}
\label{sec:stable-equivalences}

	Throughout this section, all algebras are assumed to be basic.
	We write $\underline{\mathrm{mod}}\,A$ for the category of finite-dimensional left $A$-modules modulo morphisms factoring through projective modules.
	Algebras $A$ and $B$ are \defn{stably equivalent} if $\underline{\mathrm{mod}}\,A$ and $\underline{\mathrm{mod}}\,B$ are equivalent as $\Bbbk$-linear categories.

	\begin{corollary}
		\label{cor:stable-special-multiserial}
		Let $A$ and $B$ be stably equivalent algebras.
		Then $T(A)$ is special multiserial if and only if $T(B)$ is special multiserial.
	\end{corollary}

	\begin{proof}
		By \cite[Theorem~6.4]{CX26}, $A$ is almost gentle if and only if $B$ is almost gentle.
		The assertion follows from Theorem~\ref{main:trivial-extension}.
	\end{proof}

	Recall that $A$ and $B$ are \defn{stably equivalent of Morita type} if there are bimodules ${}_AM_B$ and ${}_BN_A$, projective on both sides, with bimodule isomorphisms
	\[
		M\otimes_B N\cong A\oplus P,
		\qquad N\otimes_A M\cong B\oplus Q,
	\]
	where $P$ and $Q$ are projective bimodules over $A$ and $B$, respectively.
	The tensor functors induce mutually inverse stable equivalences.
	Such equivalences do not in general imply that $T(A)$ and $T(B)$ are stably equivalent of Morita type~\cite{LZZ17}. 

	The next corollary describes what happens when $T(A)$ is special multiserial.
	Let $s(A)$ denote the number of isomorphism classes of simple $A$-modules.

	\begin{corollary}
		\label{cor:morita-type-trivial-extensions}
		Let $A$ and $B$ be stably equivalent of Morita type, and suppose that $T(A)$ is special multiserial.
		The following are equivalent:
		\begin{enumerate}
			\item $T(A)$ and $T(B)$ are stably equivalent of Morita type;
			\item $s(A)=s(B)$;
			\item $A\cong B$.
		\end{enumerate}
		In particular, if neither $A$ nor $B$ has a semisimple direct summand, then $A\cong B$ and $T(A)\cong T(B)$, and every stable equivalence of Morita type between $A$ and $B$ lifts to a Morita equivalence.
	\end{corollary}

	\begin{proof}
		We prove $(1)\Rightarrow(2)\Rightarrow(3)\Rightarrow(1)$.
		By Theorem~\ref{main:trivial-extension}, $A$ is almost gentle, hence quadratic monomial.
		By \cite[Theorem~4.15(1)]{CX26}, we may write
		\[
			A\cong \Gamma\times\Bbbk^r,
			\qquad B\cong \Gamma\times\Bbbk^s,
		\]
		where $\Gamma$ has no semisimple direct summands and $r,s\geq0$.
		The factor $\Gamma$ is omitted if $A$ and $B$ are semisimple.

		\smallskip
		\noindent\underline{$(1)\Rightarrow(2)$}\quad
		Note that the stable centre $Z^{\mathrm{st}}$ is invariant under stable equivalence of Morita type~\cite[Proposition~5.4]{Broue94} and the dimension formula recalled in \cite[Section~3]{LZZ17} gives
		\[
			\dim_\Bbbk Z^{\mathrm{st}}\bigl(\Bbbk[\varepsilon]/(\varepsilon^2)\bigr)
			=\begin{cases}
				1,&\operatorname{char}\Bbbk\neq2,\\
				2,&\operatorname{char}\Bbbk=2.
			\end{cases}
		\]
		Since $T(\Bbbk)\cong\Bbbk[\varepsilon]/(\varepsilon^2)$, and both trivial extensions and stable centres respect finite products, we obtain
		\[
			\begin{aligned}
				0
				&=\dim_\Bbbk Z^{\mathrm{st}}(T(A))
				-\dim_\Bbbk Z^{\mathrm{st}}(T(B))=(r-s)\dim_\Bbbk Z^{\mathrm{st}}
				\bigl(\Bbbk[\varepsilon]/(\varepsilon^2)\bigr).
			\end{aligned}
		\]
		Thus $r=s$, so $s(A)=s(B)$.

		\smallskip
		\noindent\underline{$(2)\Rightarrow(3)$}\quad
		The decompositions give $s(A)-s(B)=r-s$.
		Hence $r=s$ and $A\cong B$.

		\smallskip
		\noindent\underline{$(3)\Rightarrow(1)$}\quad
		An isomorphism $A\cong B$ gives $T(A)\cong T(B)$.

		If neither algebra has a semisimple direct summand, then $r=s=0$, and the lifting assertion follows from \cite[Theorem~4.15(1)]{CX26}.
	\end{proof}

	\begin{remark}
		Note that $A$ and $A\times\Bbbk$ are stably equivalent of Morita type, whereas
		\[
			T(A\times\Bbbk)\cong T(A)\times \Bbbk[\varepsilon]/(\varepsilon^2).
		\]
		Their stable centres differ in dimension by $\dim_\Bbbk Z^{\mathrm{st}}(\Bbbk[\varepsilon]/(\varepsilon^2))>0$, so $T(A)$ and $T(A\times\Bbbk)$ are not stably equivalent of Morita type. 
	\end{remark}

	\paragraph{Use of AI.}
		The authors used GPT-6 Astra for language editing and the revision of proofs.
		In particular, GPT-6 Astra provided the proof of Lemma~\ref{lem:monomial-grading-presentation} using conjugacy of maximal tori in linear algebraic groups.
		The authors take full responsibility for the content of this article.

	\begingroup\linespread{1}\selectfont
		\bibliographystyle{alpha}
		\bibliography{reference}
	\endgroup
\end{document}